\documentclass[conference]{IEEEtran}
\usepackage{amssymb,amsthm,multicol}
\usepackage{url}
\usepackage[cmex10]{amsmath}
\ifCLASSINFOpdf
\else
\fi

\begin{document}

\newtheorem{thm}{Theorem}[section]
\newtheorem{cor}{Corollary}[section]
\newtheorem{lemm}{Lemma}[section]
\newtheorem{conj}{Conjecture}[section]
\theoremstyle{definition}
\newtheorem{defn}{Definition}[section]
\newtheorem{exmpl}{Example}[section]
\newtheorem{remrk}{Remark}[section]

%
\title{On the Components\\ 
of an Odd Perfect Number}

\author{\IEEEauthorblockN{Jose Arnaldo B. Dris}
\IEEEauthorblockA{Candidate for Master of Science in Mathematics\\
De La Salle University\\
Manila, Philippines 1004\\
Email: Jose.Dris@Safeway.Com}
}

\maketitle

\begin{abstract}
\boldmath
If $N = {p^k}{m^2}$ is an odd perfect number with special prime factor $p$, then it is proved that ${p^k} < \frac{2}{3}{m^2}$.  Numerical results on the abundancy indices $\frac{\sigma(p^k)}{p^k}$ and $\frac{\sigma(m^2)}{m^2}$, and the ratios $\frac{\sigma(p^k)}{m^2}$ and $\frac{\sigma(m^2)}{p^k}$, are used.  It is also showed that $m^2 > \frac{\sqrt{6}}{2}({10}^{150})$.
\end{abstract}


\section{INTRODUCTION}
A natural number $N$ is said to be \emph{perfect} if $N$ is equal to the sum of its proper divisors.  That is, $N$ is perfect if $\sigma(N) = 2N$.  All known perfect numbers are \emph{even}, and are generated by the formula $N = \frac{1}{2}{M_p}({M_p} + 1)$ where $M_p = {2^p} - 1$ is prime \cite{Voight}.  As of September 2006, forty-four (44) even perfect numbers have been determined, with the largest one corresponding to the Mersenne prime $2^{32582657} - 1$ \cite{GIMPS}.\\

The problem of the existence of \emph{odd perfect numbers} remains unsolved to this day.  While nobody has been able to come up with an example of an odd perfect number, no one also has been able to prove that none exists, although it is possible to derive necessary conditions for their existence.  For instance, an age-old result by Euler states that an odd perfect number, if it exists, must take the form $N = {p^k}{m^2}$ where $p \equiv k \equiv 1 \pmod 4$ and $\gcd(p, m) = 1$.  (Here, $p$ is called the special prime factor of $N$.)  The assertion that "\emph{There does not exist an odd perfect number}" has come to be called the OPN Conjecture.\\

In April of 2006, \cite{OddPerfect} was able to complete the last difficult factorization to prove that an odd perfect number $N > {10}^{500}$ using extensive computer calculations.  The other factorizations required for the proof tree are considered easy, and we may reasonably expect a proof for $N > {10}^{500}$ anytime soon.\\

There are, however, no known results even for the simple problem of comparing the size of $p^k$ to that of $m^2$.  In Dr.~ Douglas Iannucci's opinion, such a result would be difficult to obtain.  A related result due to Starni \cite{Starni} is the following:  If $N = {p^k}{m^2}$ is an odd perfect number with special prime $p$, $k + 2$ is prime, and $k + 2$ is relatively prime to $p - 1$, then $k + 2$ divides $m^2$.  A most helpful result, one which would also be difficult to obtain, would be an upper bound on $k$.
  

\section{RESULTS}
To effectively compare $p^k$ and $m^2$, we take a three-step approach:
\begin{itemize}
\item{1. Compute bounds for the abundancy indices $\frac{\sigma(p^k)}{p^k}$ and $\frac{\sigma(m^2)}{m^2}$.}
\item{2. Show that $\frac{\sigma(m^2)}{p^k}$ is an integer.}
\item{3. Determine a suitable lower bound for $\frac{\sigma(m^2)}{p^k}$.}
\end{itemize}
In Step \#3, a result found in \cite{Dandapat} is used.

\subsection{STEP \#1}
We use the fact that $p$ is a prime with $p \equiv 1 \pmod 4$ to prove the following lemma:

\begin{lemm}\label{Lemma1}
Let $N = {p^k}{m^2}$ be an odd perfect number with special prime factor $p$.  Then
$$1 < \frac{\sigma(p^k)}{p^k} < \frac{5}{4} < \frac{8}{5} < \frac{\sigma(m^2)}{m^2} < 2$$
\end{lemm}

\begin{IEEEproof}
The first inequality follows from the fact that $p$ is a prime number and $\sigma(n) \ge n, \forall n$, with equality occurring only when $n = 1$.  The second inequality we prove as follows:
$$\frac{\sigma(p^k)}{p^k} = \frac{p^{k + 1} - 1}{{p^k}(p - 1)} < \frac{p^{k + 1}}{{p^k}(p - 1)} = \frac{p}{p - 1} = \frac{1}{1 - \frac{1}{p}}$$

But since $p$ is a prime number congruent to 1 modulo 4, then $p \ge 5$.  This implies:

$$\frac{1}{1 - \frac{1}{p}} \le \frac{5}{4}$$

which means that

$$\frac{\sigma(p^k)}{p^k} < \frac{5}{4}$$

The fourth and fifth inequalities are proved by noting that

$$\frac{\sigma(p^k)}{p^k}\frac{\sigma(m^2)}{m^2} = 2$$
\end{IEEEproof}

\subsection{STEP \#2}
From the equality 
$$\frac{\sigma(p^k)}{p^k}\frac{\sigma(m^2)}{m^2} = 2$$
we consider the slightly different representation
$$\frac{\sigma(p^k)\sigma(m^2)}{p^k} = 2{m^2}$$
to arrive at the following lemma:

\begin{lemm}\label{Lemma2}
Let $N = {p^k}{m^2}$ be an odd perfect number with special prime factor $p$.  Then
$p^k$ divides $\sigma(m^2)$.
\end{lemm}

\begin{IEEEproof}
First, use the Euclidean Algorithm to determine $\gcd({p^k},\sigma(p^k))$:

$$\sigma(p^k) = \frac{p^{k + 1} - 1}{p - 1}$$
$$(p - 1)\sigma(p^k) = (p - 1){p^k} + ({p^k} - 1)$$
$$\sigma(p^k) = {p^k} + \sigma(p^{k - 1})$$
$$p^k = (p - 1)\sigma(p^{k - 1}) + 1$$
(Note that $\sigma(p^{k - 1}) < {p^k}$).  The last nonzero remainder is $\gcd({p^k},\sigma(p^k)) = 1$.  Now, $\frac{\sigma(p^k)\sigma(m^2)}{p^k} = 2{m^2}$.  Since $p^k$ is relatively prime to $\sigma(p^k)$, then $p^k$ divides $\sigma(m^2)$.
\end{IEEEproof}

\subsection{STEP \#3}
We need the following intermediate result:

\begin{lemm}\label{Lemma3}
$\sigma(A^2)$ is odd for any natural number $A$.
\end{lemm}

\begin{IEEEproof}
Let $A = \displaystyle\prod_{i = 1}^{r}{{p_i}^{\alpha_i}}$ be the prime factorization of $A$.  Then $A^2 = \displaystyle\prod_{i = 1}^{r}{{p_i}^{2{\alpha_i}}}$, and
$$\sigma(A^2) = \sigma(\displaystyle\prod_{i = 1}^{r}{{p_i}^{2{\alpha_i}}}) = \displaystyle\prod_{i = 1}^{r}{\sigma({p_i}^{2{\alpha_i}})}$$ $$= \displaystyle\prod_{i = 1}^{r}{\left(1 + p_i + {p_i}^2 + \ldots + {p_i}^{2{\alpha_i}}\right)}$$
This last product is odd regardless of whether $A$ is odd or even.
\end{IEEEproof}

Lemmas \ref{Lemma2} and \ref{Lemma3} imply that $\frac{\sigma(m^2)}{p^k}$ is odd.  In particular, $\frac{\sigma(m^2)}{p^k} \ne 2$. \\

Suppose $\frac{\sigma(m^2)}{p^k} = 1$.  Then $\frac{\sigma(p^k)}{m^2} = 2$.  This means that $\sigma(m^2) = p^k$ and $\sigma(p^k) = 2{m^2}$, or $\sigma(\sigma(m^2)) = 2{m^2}$ (i.e., $m^2$ is superperfect).  However, we have the following 1975 result from \cite{Dandapat}:

\begin{thm}\label{Theorem1}
No odd perfect number $N = {p^k}{m^2}$ satisfies $\sigma(m^2) = p^k$ and $\sigma(p^k) = 2{m^2}$.
\end{thm}

Theorem \ref{Theorem1} and the previous considerations imply that $\frac{\sigma(m^2)}{p^k} \ge 3$.  We use this lower bound in the next section to prove our main result.

\section{CONCLUSION}
We now have the following theorem:

\begin{thm}\label{Theorem2}
Let $N = {p^k}{m^2}$ be an odd perfect number with special prime factor $p$.  Then $p^k < \frac{2}{3}{m^2}$.
\end{thm}

\begin{IEEEproof}
From Lemma \ref{Lemma1}, $\sigma(m^2) < 2{m^2}.$  From Step \#3, $\sigma(m^2) \ge 3{p^k}$.  The result readily follows.
\end{IEEEproof}

Using the lower bound $N > {10}^{300}$ \cite{Brent} and Theorem \ref{Theorem2}, we arrive at the following lower bound for $m^2$:

\begin{cor}\label{Corollary1}
Let $N = {p^k}{m^2}$ be an odd perfect number with special prime factor $p$.  Then $m^2 > {\frac{\sqrt{6}}{2}}({10}^{150})$.
\end{cor}

\section{SOME NOTES}
This paper was published in the Electronic Proceedings of the 9th De La Salle University - Science and Technology Congress on July 4, 2007.  The results contained herein form part of the author's M.~Sc.~ thesis (available online via \url{http://arxiv.org/abs/1204.1450}), which was completed in August of 2008.  The author's thesis adviser (Dr.~ Gervacio) encouraged him to submit a paper containing his preliminary results to the S\&T Congress in preparation for his M.~Sc.~ thesis proposal and final defenses.  

\section*{ACKNOWLEDGMENTS}
The author would like to thank the following for their help: \\

Dr.~ Carl Pomerance of the Department of Mathematics, Dartmouth College, Hanover, NH; \\

Dr.~ Douglas Iannucci of the Science and Math Division, University of the Virgin Islands, St. Thomas, VI; \\

Dr.~ Severino Gervacio of the Mathematics Department, De La Salle University, Manila; \\

Dr.~ Blessilda Raposa of the Mathematics Department, De La Salle University, Manila; \\

Dr.~ Fidel Nemenzo of the Department of Mathematics, University of the Philippines, Quezon City; \\

Dr.~ Julius Basilla of the Department of Mathematics, University of the Philippines, Quezon City; \\

Mr. Christopher Thomas Cruz. \\

\end{document}